\documentclass[11pt,A4paper]{amsart}

\usepackage{amssymb,amsmath,eufrak}
\author[Florent Benaych-Georges]{Florent Benaych-Georges}
\title[Taylor expansions of $R$-transforms]{Taylor expansions of $R$-transforms, application to supports and  moments}
\date{\today}

\newcommand{\NC}{\operatorname{NC}}

\newcommand{\Mob}{\operatorname{Mob}}
\newcommand{\ds}{\displaystyle}

\newcommand{\ssi}{if and only if }

\newcommand{\teo}{theorem }
\newcommand{\teov}{theorem, }

\newcommand{\cids}{$\ast$-infinitely divisible distributions }
\newcommand{\fid}{$\sst\boxplus$-infinitely divisible distribution }
\newcommand{\fids}{$\sst\boxplus$-infinitely divisible distributions }

\newcommand{\R}{\mathbb{R}}
\newcommand{\C}{\mathbb{C}}

\newcommand{\z}{\mathbb{Z}}

\newcommand{\ud}{\mathrm{d}}

\newcommand{\pro}{probability }

\newcommand{\f}{\frac}
\newcommand{\ff}{\frac{1}}
\newcommand{\lf}{\left}
\newcommand{\ri}{\right}

\newcommand{\st}{such that }

\newcommand{\eps}{\varepsilon}

\newcommand{\vfi}{\varphi}
\newcommand{\ste}{\, ;\, }
\newcommand{\mc}{\mathcal }

\newcommand{\sst}{\scriptstyle}

\newcommand{\bxp}{{\scriptstyle\boxplus}}
\newcommand{\bxpp}{{\scriptscriptstyle\boxplus}}

\newcommand{\h}{\mc{H}}
\newtheorem{Th}{Theorem}[section]
\newtheorem{propo}[Th]{Proposition} 
 
\newtheorem{lem}[Th]{Lemma}

\newtheorem{rmq}[Th]{Remark}
\newtheorem{cor}[Th]{Corollary}

\newtheorem{Thappendix}{Theorem}
\newtheorem{propoappendix}[Thappendix]{Proposition} 
 
\newtheorem{lemappendix}[Thappendix]{Lemma}

\newenvironment{pr}{\noindent {\bf Proof. }}{\ \ \ $\square$}

\long\def\symbolfootnote[#1]#2{\begingroup%
\def\thefootnote{\fnsymbol{footnote}}\footnote[#1]{#2}\endgroup} 
                 
\begin{document}


\maketitle
\symbolfootnote[0]{{\it MSC 2000 subject classifications:}  primary 60E10, 
46L54, 
secondary 60E07} 
 
\symbolfootnote[0]{{\it Key words:} $R$-transform, free cumulants, free infinitely divisible distributions} 

\begin{abstract}We prove that a probability measure on the real line has a moment of order $p$ (even integer), \ssi its $R$-transform admits a Taylor expansion with $p$ terms. We also prove a weaker version of this result when $p$ is odd. We then apply this to prove that a \pro measure whose $R$-transform extends analytically to a ball with center zero is compactly supported, and that a free infinitely divisible distribution has  a moment of order $p$ even, \ssi its L\'evy measure does so. We also prove a weaker version of the last result when $p$ is odd.\end{abstract}

\section*{Introduction} 
Addition of free random variables gives rise to a convolution $\bxp$ on the set of \pro measures on the real line. The operation $\bxp$, defined in \cite{defconv}, is called the free convolution. The classical convolution $*$ is linearized by the logarithm of the Fourier transform: the Fourier transform of the classical convolution of \pro measures is the product of their Fourier transforms. In the same way, the $R$-transform of the free convolution of \pro measures is the sum of their $R$-transforms. 
The existence of moments of even order of a \pro measure is linked to the Taylor expansion of its Fourier transform in the neighborhood of zero (\cite{feller}, XV.9.15). Moreover, the coefficients of the logarithm of this Taylor expansion are, up to a division by a factorial, the classical cumulants of the measure (\cite{mattner}). In this paper, we prove  that a probability measure has a moment of even order $p$ \ssi its $R$-transform admits a Taylor expansion with $p$ terms. Moreover, in this case, the coefficents of this expansion are the first $p$ free cumulants of the measure (defined in \cite{spei}). When $p$ is odd, one implication (moment $\Rightarrow$ Taylor expansion) stays true, and the other one is maintained when the support of the measure is minorized or majorized. These results can be transferred to the Voiculescu transform. 

The first consequence of this result is a criterion for compactness of the support of a measure. Roughly speaking (details in section \ref{riderzpremierjour}), the $R$-transform of a \pro measure is the analytic inverse of its Cauchy transform on the intersection of a cone with origin zero and a ball with center zero.   For compactly supported measures, the inversions can be done on balls centered at zero without intersecting them with cones. But it had not yet been proved that any $R$-transform which, once defined on the intersection of a cone and a ball,  can be analytically extended to the whole ball, is the one of a compactly supported \pro measure. We prove it here, and our result about Taylor expansions of $R$-transforms seems necessary to prove it.

We also apply this result to prove that a \fid has  a moment of order $p$ even \ssi its L\'evy measure does so. As before, when $p$ is odd, one implication is maintained ($p^{th}$ moment for L\'evy mesure $\Rightarrow$ $p^{th}$ moment for the distribution), and the other one is maintained under the additional assumption that the support of the L\'evy measure is minorized or majorized.   As this paper was already written, Thierry Cabanal-Duvillard (\cite{cabduv.ID}), using random matrices,  proved another result of this type. 

The first section of the paper is devoted for the presentation of the tools, for the proof of the result about Taylor expansions of $R$-transforms (\teo \ref{13.08.04.1}) and for the criterion that follows. The proof of the theorem relies on the Hamburger-Nevanlinna \teo and on proposition \ref{invdeDLdansH} of the appendix. In the second part, we apply the result to $\bxp$-infinitely divisible distributions.

{\bf Aknowledgements.} The author would like to thank his advisor Philippe Biane, as well as Professor Hari Bercovici for his encouragements.  Also, he would like to thank  C\'ecile Martineau for her contribution to the English version of this paper. 

\section{Asymptotic expansions of the $R$-transform}\label{riderzpremierjour}
Let us first present the $R$-transform (for further details, see \cite{defconv}). We define, for $\alpha, \beta >0$, the set \begin{eqnarray*}
\Delta_{\alpha,\beta}&=&\{z=x+iy\ste |x|<-\alpha y, |z|<\beta\}.
\end{eqnarray*} 
Note that $z\mapsto 1/\!z$ maps the set denoted by $\Gamma_{\alpha,\ff{\beta}}$ in \cite{defconv},\cite{appenice} onto $\Delta_{\alpha,\beta}$. In order to make notations lighter, we have prefered to avoid the references to the sets $\Gamma_{.,.}$, hence to compose the Cauchy transform on the right with the map $z\mapsto 1/\!z$.   
                                                                                                                                                                                                                                                      
The {\it Cauchy tranform} of a probability measure $\mu$ on the real line is the analytic function on the upper half-plane $G_\mu : z \mapsto \int_\R\f{\ud \mu(t)}{z-t}$. 

By Proposition 5.4 and Corollary 5.5 of \cite{defconv}, for all positive numbers $\alpha$, for all $\eps \in (0,\min\{\alpha,1\})$, for $\beta$ small enough, $z\mapsto G_\mu(1/\!z)$ is a conformal bijection from $\Delta_{\alpha,\beta}$ onto an open set $\mc{D}_{\alpha,\beta}$, \st $$\Delta_{\alpha-\eps,(1-\eps)\beta }\subset \mc{D}_{\alpha,\beta}  \subset  \Delta_{\alpha+\eps,(1+\eps)\beta}.$$ 

The inverses of $G_\mu(1/\!z)$ on all $   \mc{D}_{\alpha,\beta}$'s define together an analytic function $L_\mu$ on the union $\mc{D}$ of their domains. This function is a right inverse of $G_\mu(1/\!z)$ on $\mc{D}$. Moreover, the open set $\{z\ste G_\mu(1/\!z)\in \mc{D}\}$ has a unique connected component which contains a set of the type $\Delta_{\alpha,\beta}$, and on this connected component, $G_\mu(1/\! z)$ is also a right inverse of $L_\mu$ (it is true on a set of the type $\Delta_{\alpha,\beta}$ by Proposition 5.4 of \cite{defconv}, and therefore it is true on the connected component by analycity).  

The {\it $R$-transform} of $\mu$ is $R_\mu(z)=K_\mu(z)-1/\!z$, where $K_\mu=1/\!L_\mu$.   
 
The natural space for $R$-transforms is the space, denoted by $\h$, of functions $f$ which are analytic in a domain $\mc{D}_f$ \st for all positive $\alpha$, there exists a positive $\beta$ \st $$\Delta_{\alpha,\beta}\subset \mc{D}_f.$$ The introduction of this space is not necessary to work with $R$-transforms, but it will be usefull in our work  on their Taylor expansions.

One can summarize the different steps of the construction of the $R$-transform in the following chain$$
\begin{array}{l}\ds\underset{\substack{{\textrm{probability}}\\ {\textrm{measure}}}}{\mu}\,\,\longrightarrow \,\,\underset{\substack{{\textrm{Cauchy}}\\ {\textrm{transf.}}}}{G_\mu}\,\,\longrightarrow\,\, \underset{\textrm{function of $\h$}}{L_\mu(z)=\lf(G_\mu(1/\!z)\ri)^{-1}}\,\,\longrightarrow\,\, \underset{\substack{{\textrm{function}}\\ {\textrm{of $\h$}}}}{K_\mu=1/\!L_\mu}\,\, \longrightarrow\\ \\
\ds\underset{\textrm{function of $\h$}}{R_\mu(z)=K_\mu(z)-1/\!z.}\end{array}$$

For example, for respectively $\mu=\delta_a$, $\f{2}{\pi r^2}\sqrt{r^2-(x-m)^2}1_{|x-m|\leq r}\ud x$, $\f{\ud x}{\pi(1+x^2)}$, $R_\mu(z)=a$, $m+\f{r^2}{4}z$, $-i$.
  
  The main property of the $R$-transform is the fact that it linearizes free convolution: for all $\mu,\nu$, $R_{\mu\bxpp\nu}=R_\mu+R_\nu$. It is also useful for characterizing tight sets of \pro measures and for giving a necessary and sufficient condition for weak convergence (\cite{defconv}). 
\begin{rmq} Some authors prefer to work with the {\it Voiculescu transform} $\vfi_\mu(z)=R_\mu(1/\!z)$ rather than with the $R$-transform. Using the fact that $z\mapsto 1/\!z$ maps the set denoted  by $\Gamma_{\alpha,\ff{\beta}}$ in \cite{defconv}, \cite{appenice}  onto $\Delta_{\alpha,\beta}$, all our results can be transferred to Voiculescu transforms.\end{rmq}

One can wonder how to express the $R$-transform of a measure more directly from the measure. In the case where $\mu$ is compactly supported, one can extend its Cauchy transform to the complementary of a closed ball with center zero, and repeat the previous work  replacing every  $\Delta_{\alpha,\beta}$ by the ball $\bar{B}(0,\beta)$. It follows that the $R$-transform can be defined analytically to a neighborhood of zero, and it has been proved in \cite{spei} that the coefficients of its series expansion \begin{equation}\label{seriesexpRtr}R_\mu(z)=\sum_{i\geq 0}k_{i+1}(\mu)z^i\end{equation}are the {\it free cumulants} of $\mu$, defined by any of the two equivalent formulas:
\begin{eqnarray}\ds \forall i\geq 1,\quad m_i(\mu)&=&\sum_{\pi\in \NC(i)}\;\prod_{\substack{
\textrm{$V$ class}\\ 
\textrm{of $\pi$}}}
k_{|V|}(\mu),\label{cum->mom}\\ 
\ds \forall i\geq 1,\quad k_i(\mu)&=&\sum_{\pi\in \NC(i)}\;\Mob(\pi)\prod_{\substack{
\textrm{$V$ class}\\ 
\textrm{of $\pi$}}}
m_{|V|}(\mu),\label{defcum}\end{eqnarray}
where for all integers $i$, $m_i(\mu)$ is the $i$-th moment of $\mu$,   $\NC(i)$ denotes the set of partitions of $\{1$,..., $i\}$ \st there does not exist $1\leq j_1<j_2<j_3<j_4\leq i$ with $j_1,j_3$ in the same class and $j_2,j_4$ in another class (such partitions are said to be {\it non crossing}), and  $\Mob$ is a function on $\NC(i)$ which we will not need to explicit here, called the M\"obius function (for further details, see \cite{spei}).

\begin{rmq}The {\it classical cumulants} $c_i(\mu)$ of a compactly supported probability measure $\mu$ can be defined by the analogous equation (\cite{mattner}):  
$$\ds \forall i\geq 1,\quad m_i(\mu)=\sum_{
\substack{
\textrm{$\pi$ partition}\\ 
\textrm{of $\{1$,..., $i\}$}
}}\;
\prod_{
\substack{
\textrm{$V$ class}\\ \textrm{of $\pi$}
}}c_{|V|}(\mu),$$ and one has, for all $z$ complex numbers, $$\ds \exp \sum_{i\geq 1}\f{c_i(\mu)}{i!}z^i=\int_\R e^{tz}\ud \mu(t).$$\end{rmq}

If $\mu$ is non compactly supported, but admits a $p^{th}$ moment (hence, by the H\"older inequality, moments of order $1$,..., $p$), one can define its  first $p$ cumulants by the equation (\ref{defcum}) for $i=1$,..., $p$. Theorem \ref{13.08.04.1} bellow extends  (\ref{seriesexpRtr}) to this case.

First of all, for functions of $\h$, we only consider  ``non tangential limits'' to zero. That is, if $f\in\h$, $$\ds\lim_{z\to 0} f(z)=l$$ means that  for all positive $\alpha$, for a certain $\beta$,  $$\ds\lim_{\substack{z\to 0\\z\in \Delta_{\alpha,\beta}}} f(z)=l,$$ and {\it $f$ admits a Taylor expansion of order $p$} means that there exist complex numbers $a_0,\ldots, a_p$, a function $\upsilon\in \h$, \st $$\ds f(z)=\sum_{i=0}^p a_iz^i+z^p\upsilon(z), \textrm{ with }\lim_{z\to 0}\upsilon(z)=0.$$ A subset of $\R$ will be said to be {\it minorized} (resp. {\it majorized}) if it is contained in an interval of the type $(a,+\infty)$ (resp. of the type $(-\infty,a)$),  where $a$ is a real number.  
\begin{Th}\label{13.08.04.1}Let $p$ be a positive  integer and $\mu$ be a \pro measure on the real line.\\
(a) If $\mu$ admits a $p^{th}$ moment, then $R_\mu$ admits the Taylor expansion 
$$\ds R_\mu(z)=\sum_{i=0}^{p-1}k_{i+1}(\mu)z^i+o(z^{p-1}).$$
(b) Conversely, if  $p$ is even or if the support of $\mu$  is minorized or majorized, and if $R_\mu$ admits a Taylor expansion of order $p-1$ with real coefficients, then $\mu$ has a $p^{th}$ moment.
\end{Th}

\begin{rmq} In the second part of the \teov the coefficents have to be real. For example, the Cauchy distribution $\mu=\f{\ud x}{\pi(1+x^2)}$, has no moments, while $R_\mu=-i$.
\end{rmq}

The proof of the \teo uses the work on functions of $\h$ done in the appendix, and the Hamburger-Nevanlinna \teov that we give here:
\begin{Th}\label{13.08.04.2}Let $p$ be a positive integer and $\mu$ be a \pro measure on the real line.\\
(i) If $\mu$ admits a $p^{th}$ moment, then $G_\mu(1/\!z)$ admits the Taylor expansion
$$\ds G_\mu(1/\!z)=\sum_{i=1}^{p+1}m_{i-1}(\mu)z^{i}+o(z^{p+1}).$$
(ii) Conversely, if $p$ is even or if the support of $\mu$  is minorized or majorized, and if $G_\mu(1/\!z)$ admits the Taylor expansion of order $p+1$ with real coefficients, then $\mu$ has a $p^{th}$ moment.
\end{Th}

The part $(i)$ and the part $(ii)$ in the case where $p$ is even are respectevely the first and second parts of Theorem 3.2.1 of 
\cite{akhi}. There is no reference for the part $(ii)$ in the case where  the support of $\mu$  is minorized or majorized, so we prove it:

\begin{pr} We suppose that the  support of $\mu$  is minorized or majorized, and that there exists real numbers $r_0,\ldots, r_{p+1}$ \st $G_\mu(1/\!z)$ admits the Taylor expansion
$$\ds G_\mu(1/\!z)=\sum_{i=0}^{p+1}r_{i}z^{i}+o(z^{p+1}).$$
First of all, by Proposition 5.1 of \cite{defconv}, $r_0=0$ and $r_1=1$. Let us prove, 
by induction over $q\in \{0,\ldots, p\}$, that $\mu$ admits a $q^{th}$ moment and that for all 
$l=0,\ldots, q$, $m_l(\mu)=r_{l+1}$. For $q=0$, it is obvious by what precedes. Suppose that the result has been proved to  rank $q-1$.  We have 
$$r_{q+1}
=\ds\lim_{\substack{y\to 0\\ y>0}}\ff{(iy)^{q+1}}\lf(G_\mu\lf(\ff{iy}\ri)-\sum_{l=1}^{q}r_l(iy)^l\ri).$$Using the fact that for all $l=1,\ldots,q$, $r_l=m_{l-1}$, and taking the real part, we obtain 
$$r_{q+1}=\ds\lim_{\substack{y\to 0\\ y>0}}\int \f{t^{q}}{1+y^2t^2}\ud \mu(t).$$
Let us write $\mu=\mu^++\mu^-$, where $\mu^+,\mu^-$ are finite positive measures with supports respectively contained in the non negative and in the non positive real half lines. Now we have  $$r_{q+1}=\ds\lim_{y\to 0}\int \f{t^{q}}{1+y^2t^2}\ud\mu^+(t)+\int \f{t^{q}}{1+y^2t^2}\ud \mu^-(t).$$By hypothesis on $\mu$, one of the measures $\mu^+$, $\mu^-$ has compact support. We apply the theorem of dominated convergence to the term corresponding to this measure, and then the \teo of monotone convergence to the other term. The desired result follows.
\end{pr}

{\bf Proof of theorem \ref{13.08.04.1}. } We use the abbreviation {\it T.e.r.c.} for ``Taylor expansion with real coefficients''.\\ 
(a) By the Hamburger-Nevanlinna \teov $G_\mu(1/\!z)$ admits an T.e.r.c. of order $p+1$  with leading term $z$. So, by proposition \ref{invdeDLdansH} of the appendix, its inverse $1/\!K_\mu$ does also. Thus, dividing by $z$, $1/\!(zK_\mu(z))$ admits an T.e.r.c. of order $p$ with leading term $1$. 
It is obvious that if one composes on the left a function of $\h$ which admits a T.e.r.c.   of order $p$ with leading term $1$ by a function which admits a T.e.r.c.  of order $p$ in a neighborhood of $1$, one obtains a function of $\h$ which admits a T.e.r.c. of order $p$.   Therefore, $zK_\mu(z)$ does so: there exists $a_1$,..., $a_p\in \R$ \st  $$\ds zK_\mu(z)= 1+\sum_{l=1}^p a_lz^l+o(z^{p}).$$ 
So, since $R_\mu(z)=K_\mu(z)-1/\!z$, it suffices to prove that the coefficients $a_1,\ldots, a_p$ are the $p$ first free cumulants of $\mu$. 
\\
Let us denote these coefficients by $a_1(\mu), \ldots, a_p(\mu)$ (the functions $a_1,\ldots,a_p$ are defined on the set of \pro measures with $p^{th}$ moment). 
\\
It is easy to see that there exists polynomials $P_1,\ldots,P_p\in \R[X_1,\ldots,X_p]$ defined by: for all $x_1$, $y_1$,..., $x_p,y_p\in\R$, for all bijective function $f$  defined in a neighborhood of zero  \st $f(z)=z+\sum_{i=2}^{p+1}x_{i-1}z^{i}+o(z^{p+1})$ and $f^{-1}(z)= z+\sum_{i=2}^{p+1}y_{i-1}z^{i}+o(z^{p+1})$, one has for all $ i=1$,...,$p$, $y_i=P_i(x_1,\ldots,x_p)$ (the existence of these polynomials easily follows from the equation $f\circ f^{-1}(z)=z$).
\\
Thus, $$1/\!K_\mu(z)=z+\sum_{i=2}^{p+1}P_{i-1}(m_1(\mu),\ldots,m_p(\mu))z^{i}+o(z^{p+1}),$$ and the polynomials $P_1$,\ldots,$P_p$ do not depend on the choice of the \pro measure with $p^{th}$ moment $\mu$.         Continuing with the same kind of argument (where the fact that the leading term of the T.e.r.c. of $1/\!(zK_\mu(z))$ is $1$ is crucial), one sees that the functions $a_1$,..., $a_p$ can be expressed as universal polynomials in the moments: there exist $Q_1$,..., $Q_p\in \R[X_1,\ldots,X_p]$ (not depending on $\mu$) \st for all $i=1$,...,$p$, $a_i(\mu)=Q_i(m_1(\mu),\ldots, m_p(\mu))$. 

Moreover, by (\ref{defcum}), there also exists universal polynomials $R_1$,..., $R_p\in \R[X_1,\ldots,X_p]$ \st for all $i=1$,...,$p$, $k_i(\mu)=R_i(m_1(\mu),\ldots, m_p(\mu))$ (for further details on these polynomials, see \cite{spei}). 

It only remains to prove that for all $i$, $Q_i$ and $R_i$ coincide on the $p$-tuple $(m_1(\mu)$,..., $m_p(\mu))$. There are several ways to see it. First, the Theorem of the third part of \cite{spei} asserts that $Q_i=R_i$. But without refering directly to this result, by  (\ref{seriesexpRtr}), $Q_i$ and $R_i$ coincide on the $p$-tuple $(m_1(\mu)$,..., $m_p(\mu))$ when $\mu$ is compactly supported, and it is proved in the chapter of \cite{akhi} devoted to the solution of the moment problem that there exists a compactly supported \pro measure (more precisely a convex combination of Dirac measures) with $p$ first moments  $(m_1(\mu)$,..., $m_p(\mu))$. 

(b) If $R_\mu$ admits an T.e.r.c. of order $p$, we prove, with the inverse operations of the ones of the first part of the proof of (a) of this \teov that (ii) of  the Hamburger-Nevanlinna \teo applies. \ \ \ $\square$

Note that we do not know yet if any  $R$-transform which, once defined as an element of $\h$, 
can be analytically extended to a whole ball with center zero, comes from a compactly supported \pro measure. The affirmative answer will be given the following corollary.

\begin{cor}\label{ebaypremierjour} A \pro measure on the real line is compactly supported \ssi its $R$-transforms can be analytically extended to an open ball with center zero.\end{cor}
 
\begin{pr}As mentioned above, it is already known that the $R$-transform of a compactly supported \pro measure extends analytically to an open ball with center zero. Suppose now that the $R$-transform of a \pro measure $\mu$ extends analytically to an open ball with center zero and positive radius $r$. By the previous \teov $\mu$ has moments of all orders, and its free cumulants are the coefficients of the series expansion of its $R$-transform. Hence the sequence $|k_n(\mu)|^{\ff{n}}$ is bounded by a number $C$, and so, from (\ref{cum->mom}) and the majorations  $\#\NC(n)\leq 4^n$,   $\forall \pi\in \NC(n)$, $|\Mob(\pi)|\leq 4^n$ (\cite{pfi}, p. 149-150),  one has, for all integers $n$, $|m_n(\mu)|\leq (16C)^n$. This implies that the support of $\mu$ is contained in $[-16C,16C]$: otherwise, there exists $\eps,\delta>0$ \st $\mu(\R-[-16C-\delta, 16C+\delta])>\eps$, and so $\forall n, m_{2n}(\mu)>\eps(16C+\delta)^{2n}$, which contradicts  $|m_{2n}(\mu)|\leq (16C)^{2n}$. \end{pr}

\section{Moments of $\bxp$-infinitely divisible distributions}
First, recall that a \pro measure $\mu$ is said to be $\bxp$-infinitely divisible (resp.  $\ast$-infinitely divisible) if for all integer $n$, there exists a \pro measure $\nu_n$ \st $\nu_n^{\bxpp n}=\mu$ (resp.   $\nu_n^{\ast n}=\mu$). This condition is equivalent (\cite{appenice}) to the existence of a sequence $(\mu_n)$ of \pro measures \st $\mu_n^{\bxpp n}$ (resp.  $\mu_n^{\ast n}$) converges weakly to $\mu$. These distributions have been classified in \cite{defconv} (resp. in \cite{gne}): $\mu$ is $\bxpp$ (resp. $*$)-infinitely divisible \ssi there exists a real number $\gamma$ and a positive finite measure on the real line (abbreviated from now on into {\it p.f.m.}) $\sigma$ \st $R_\mu(z)=\gamma+\int_\R\f{z+t}{1-tz}\ud \sigma(t)$ (resp. the Fourier transform is $\hat{\mu}(t)=\exp\lf[i\gamma t+\int_\R(e^{itx}-1-\f{itx}{x^2+1})\f{x^2+1}{x^2}\ud \sigma(x)\ri]$). Moreover, in this case, such a pair $(\gamma, \sigma)$ is unique, and we denote $\mu$ by $\nu_{\bxpp\!}^{\gamma, \sigma}$ (resp. $\nu_{\ast}^{\gamma, \sigma}$). Thus, one can define a bijection from the set of \cids to the set of \fids by $\nu_{\ast}^{\gamma, \sigma}\mapsto \nu_{\bxpp\!}^{\gamma, \sigma}$. This bijection, called the Bercovici-Pata bijection, is defined in a formal way, but appears to have deep properties. First of all, it is easy to see that for all $(\gamma, \sigma)$ and $(\gamma', \sigma')$, 
$$\nu_{\bxpp}^{\gamma, \sigma}\bxpp \nu_{\bxpp}^{\gamma', \sigma'}=\nu_{\bxpp\!}^{\gamma+\gamma', \sigma+\sigma'},\qquad \nu_{\ast}^{\gamma, \sigma}\ast \nu_{\ast}^{\gamma', \sigma'}=\nu_{\ast}^{\gamma+\gamma', \sigma+\sigma'}.$$ Thus, the bijection previously defined is a semi-goup morphism. Moreover, it has been proved in  \cite{steen2} that it is an homeomorphism with respect to weak convergence topology. At last, a surprising property of the Bercovici-Pata bijection was proved in \cite{appenice}: for all  sequences $(\mu_n)$ of \pro measures, the sequence  $\mu_n^{\ast n}$ tends weakly to a measure $\nu_{\ast}^{\gamma, \sigma}$ \ssi the sequence $\mu_n^{\bxp n}$ tends weakly to $\nu_{\bxpp\!}^{\gamma, \sigma}$. Note that this property does not follows from the previous ones, because the measures $\mu_n$ are not supposed to be  $*$-infinitely divisible, so one cannot apply the bijection to   $\mu_n^{\ast n}$. A somewhat more concrete realization of this surprising bijection, using random matrices, can be found in \cite{cabduv.ID} and \cite{fbg.ID}.  
 
The measure $\sigma$ is said to be the {\it L\'evy measure} of $\nu_{\ast}^{\gamma, \sigma}$ and $\nu_{\bxpp\!}^{\gamma, \sigma}$. It is well known (section 25 of \cite{sato}) that a $\ast$-infinitely divisible distribution admits moments of the same orders as its L\'evy measure. We will prove, in this section, that in the free case, if the L\'evy measure has a moment of order $p$, then so does the distribution, and that the converse is true when considering moments of even order or L\'evy measures with minorized or majorized support. In a recent preprint (\cite{cabduv.ID}), Thierry Cabanal-Duvillard has proved  the first implication.  

First of all, note that a $\bxp$-infinitely divisible distribution has compact support \ssi its L\'evy measure does so. It was never written like this, but it was proved (\cite{BV92}, \cite{hiai}) that a compactly supported \pro measure $\mu$  is $\bxp$-infinitely divisible \ssi its $R$-transform can be written $\gamma+\int_\R\f{z+t}{1-tz}\ud \sigma(t)$, with $\gamma\in \R$ and $\sigma$ a compactly supported p.f.m.. Moreover, in this case, with the series expansion of $\int_\R\f{z+t}{1-tz}\ud \sigma(t)$ and (\ref{seriesexpRtr}), it is easy to see that for all positive integer $p$, the $p^{th}$ free cumulant of $\mu$ is \begin{equation}\label{16.08.04.2}k_p(\nu_{\bxpp\!}^{\gamma,\sigma})=m_{p-2}(\sigma)+m_p(\sigma),\end{equation} where $m_{-1}(\sigma):=\gamma$.

\begin{rmq}Note that in this case,  $m_{p-2}(\sigma)+m_p(\sigma)$ is also the $p^{th}$ classical cumulant of the $\ast$-infinitely divisible correspondant of $\mu$.
\end{rmq}

In the proof of proposition \ref{moments:loi<->Levy.15.08.04}, we will need the following lemma: 
\begin{lem}\label{an40plusbifluore} If the support of an f.p.m. $\sigma$ is contained in $(0,\infty)$ (resp. in $(-\infty, 0)$), then $\nu_{\bxpp\!}^{0, \sigma}$ is concentrated on $[0,\infty)$ (resp. $(-\infty,0]$).
\end{lem}

\begin{pr} The classical version of this result is well known : if the support of  $\sigma$ is contained in $(0,\infty)$ (resp. in $(-\infty, 0)$), then $\nu_{\ast}^{0, \sigma}$ is concentrated on $[0,\infty)$ (resp. $(-\infty,0]$). Thus $\nu_{\ast}^{0, \f{\sigma}{n}}$ is concentrated on $[0,\infty)$ (resp. $(-\infty,0]$), and the same holds for  $\lf(\nu_{\ast}^{0, \f{\sigma}{n}}\ri)^{\bxp n}$. But for all $n$, $\lf(\nu_{\ast}^{0, \f{\sigma}{n}}\ri)^{\ast n}=\nu_{\ast}^{0, \sigma}$. Thus  $\lf(\nu_{\ast}^{0, \f{\sigma}{n}}\ri)^{\bxp n}$ converges weakly to $\nu_{\bxpp\!}^{0, \sigma}$, which is hence concentrated on $[0,\infty)$ (resp. $(-\infty,0]$).
\end{pr}

\begin{propo}\label{moments:loi<->Levy.15.08.04}Let $\gamma$ be a real number, $\sigma$ be an f.p.m., and $p$ be a positive integer.\\
$1^\circ$) If $\sigma$ admits a moment of order $p$, then the same holds for $\nu_{\bxpp\!}^{\gamma, \sigma}$.\\
$2^\circ$) Suppose moreover that $p$ is even or that the support of $\sigma$ is minorized or majorized. In this case, if  $\nu_{\bxpp\!}^{\gamma, \sigma}$  admits a moment of order $p$, then the same holds for $\sigma$.\end{propo}

Before beginning the proof, let us recall a few results about weak convergence of \pro measures and of p.f.m.. First, for any sequences $(\gamma_n)$ of real numbers and $(\sigma_n)$ of p.f.m., the sequence $(\nu_{\bxpp\!}^{\gamma_n, \sigma_n})$ converges weakly to a $\bxp$-infinitely divisible measure $\nu_{\bxpp\!}^{\gamma, \sigma}$ \ssi $\gamma_n$ tends to $\gamma$ and $\sigma_n$ tends weakly to $\sigma$ (\cite{steen2}). Recall that a sequence $\rho_n$ of p.f.m. (including \pro measures) converges weakly to a p.f.m. if for all continuous bounded function $f$, $\int f\ud \rho_n$ tends to $\int f\ud \rho$. In this case, combining Theorems 5.1 and 5.3  of \cite{billingsley}, one has, for all nonnegative continuous function $f$, \begin{equation}\label{duredurite} \int f\ud \rho\leq \liminf \int f \ud \rho_n. \end{equation}

{\bf Proof of proposition \ref{moments:loi<->Levy.15.08.04}. } $1^\circ$) Suppose $\sigma$ to admit a moment of order $p$. Define three f.p.m. $\sigma^-,\sigma^c,\sigma^+$ by $ \sigma^-(A)=\sigma(A\cap(-\infty,-1))$, $\sigma^c(A)=\sigma(A\cap[-1,1])$, $\sigma^+(A)=\sigma(A\cap(1,\infty))$ for all Borel set $A$. Then $\nu_{\bxpp\!}^{\gamma, \sigma}=\nu_{\bxpp\!}^{0, \sigma^-}\bxpp\nu_{\bxpp\!}^{\gamma, \sigma^c}\bxpp\nu_{\bxpp\!}^{0, \sigma^+}$. So, by Minkowski inequality in $W^*$-\pro spaces (equation (26) of \cite{Nelson}), it suffices to prove that each of $\nu_{\bxpp\!}^{0, \sigma^-}$, $\nu_{\bxpp\!}^{\gamma, \sigma^c}$, $\nu_{\bxpp\!}^{0, \sigma^+}$ admits a $p^{th}$ moment. As explained above,  $\nu_{\bxpp\!}^{\gamma, \sigma^c}$ is compactly supported, so we have reduced the problem to the case where $\gamma=0$ and the support of $\sigma$ is contained in $(-\infty,0)$ or in $(0,+\infty)$. Both cases are treated in the same way, suppose for example the support of $\sigma$ to be contained in $(0,\infty)$. Let us then define, for $n$ positive integer, the f.p.m. $\sigma_n$  by $\sigma_n(A)=\sigma(A\cap (0,n))$. By dominated convergence, $\sigma_n$ tends weakly to $\sigma$. So, by what precedes, $ \nu_{\bxpp\!}^{0, \sigma_n}$ tends weakly to  $ \nu_{\bxpp\!}^{0, \sigma}$, and by (\ref{duredurite}), $\int |x|^p \ud  \nu_{\bxpp\!}^{0, \sigma}(x)\leq \liminf \int |x|^p \ud  \nu_{\bxpp\!}^{0, \sigma_n}(x)$. But by the previous lemma,   all $\nu_{\bxpp\!}^{0, \sigma_n}(x)$ are concentrated on $[0,\infty)$, so one has  $\int |x|^p \ud  \nu_{\bxpp\!}^{0, \sigma}(x)\leq \liminf m_p(\nu_{\bxpp\!}^{0, \sigma_n}(x))$: it suffices to prove the boundness of the sequence  $(m_p(\nu_{\bxpp\!}^{0, \sigma_n}))_n$. But by (\ref{16.08.04.2}), for all $n$, for all integer $q$, the $q^{th}$ free cumulant of $\nu_{\bxpp\!}^{0, \sigma_n}$ is $m_{q-2}(\sigma_n)+m_q(\sigma_n)$ (with $m_{-1}(\sigma_n)=0$). Thus one has, for all $n$, \begin{eqnarray*}\ds m_p(\nu_{\bxpp\!}^{0, \sigma_n})&=&\sum_{\pi\in \NC(p)}\prod_{V\in \pi}k_{|V|}(\nu_{\bxpp\!}^{0, \sigma_n})\\ \ds&=&\sum_{\pi\in \NC(p)}\prod_{V\in \pi}(m_{|V|-2}(\sigma_n)+m_{|V|}(\sigma_n))\\ \ds&\leq & \sum_{\pi\in \NC(p)}\prod_{V\in \pi}(m_{|V|-2}(\sigma)+m_{|V|}(\sigma)),\end{eqnarray*}and the result is proved.
 \\
$2^\circ$)  First of all, $\nu_{\bxpp\!}^{\gamma, \sigma}\ast\delta_{-\gamma}=\nu_{\bxpp\!}^{0, \sigma}$. So one can suppose that $\gamma=0$.

Moreover, let us prove that we can replace the hypothesis ``the support of $\sigma$ is minorized or majorized'' by ``the support of $\sigma$ is contained in $(1,\infty)$''. So suppose the support of $\sigma$ to be minorized or majorized, and $\nu_{\bxpp\!}^{0, \sigma}$  to admit a moment of order $p$. Since a push-forward of $\sigma$ by $x\mapsto -x$ transforms $\nu_{\bxpp\!}^{0, \sigma}$ into its push-forward by $x\mapsto -x$, one can suppose the support of $\sigma$ to be contained in an interval $(m,\infty)$ with $m\in \R$.  Define two f.p.m. $\sigma_c,\sigma^+$
 by $\sigma_c(A)=\sigma(A\cap (m,|m|+1])$, $\sigma^+(A)=\sigma(A\cap(|m|+1,\infty))$. Then $ \nu_{\bxpp\!}^{0, \sigma}=\nu_{\bxpp\!}^{0, \sigma^c}\bxpp\nu_{\bxpp\!}^{0, \sigma^+}$, thus, by \cite{defconv}, $ \nu_{\bxpp\!}^{0, \sigma}$ is the distribution of the sum $X+Y$ of two free selfadjoint operators $X,Y$ affiliated to  a $W^*$-\pro space $(\mc{A},\vfi)$, respectively distributed according to $\nu_{\bxpp\!}^{0, \sigma^c}$, $\nu_{\bxpp\!}^{0, \sigma^+}$. By hypothesis $X+Y\in L^p(\mc{A},\vfi)$, and by compactness of the support of its distribution,   $X$ is bounded. So, by the Minkowsky inequality (\cite{Nelson}), $Y=(X+Y)-X\in L^p(\mc{A},\vfi)$. So  $\nu_{\bxpp\!}^{0, \sigma^+}$ admits a $p^{th}$ moment, and it suffices to prove that $\sigma^+$ admits a $p^{th}$ moment (because $\sigma^c$ has compact support).

So let us suppose $p$ to be even or the support of $\sigma$ to be contained in $(1,\infty)$. Suppose now that  $\nu_{\bxpp\!}^{0, \sigma}$  admits a moment of order $p$. By the first part of \teo  \ref{13.08.04.1}, $R_{\nu_{\bxpp\!}^{0, \sigma}}$ admits a Taylor expansion of order $p$, and so for all positive integer $n$, the same holds for \begin{equation}\label{bamakeetsavache15.08.04}R_{\nu_{\bxpp\!}^{0, \f{\sigma}{n}}}=\ff{n}R_{\nu_{\bxpp\!}^{0, \sigma}}.\end{equation} Then  lemma \ref{an40plusbifluore} allows us to apply the second part of \teo  \ref{13.08.04.1}, and to deduce that for all $n$, $\nu_{\bxpp\!}^{0, \f{\sigma}{n}}$ admits a $p^{th}$ moment. Moreover, the coefficients of the Taylor expansions of the $R$-transforms of $ \nu_{\bxpp\!}^{0, \sigma}$, $\nu_{\bxpp\!}^{0, \f{\sigma}{n}}$ are their free cumulants, so from  (\ref{bamakeetsavache15.08.04}), we have $$\forall i=1,\ldots, p,\qquad k_i(\nu_{\bxpp\!}^{0, \f{\sigma}{n}})=\ff{n}k_i(\nu_{\bxpp\!}^{0, \sigma}).$$ 
But it was proved in \cite{defconv} (Theorem 5.10, (iii)) that $\sigma$ is the weak limit of the sequence of f.p.m. $\lf(\f{nx^2}{1+x^2}\ud \nu_{\bxpp\!}^{0,\sigma_n}\ri)$. So, by (\ref{duredurite}), \begin{eqnarray*}\int |x|^p\ud\sigma (x)&\leq&\liminf \int\f{|x|^pnx^2}{1+x^2}\ud\nu_{\bxpp\!}^{0,\sigma_n}(x) \\
&\leq& \liminf \int|x|^pn\ud\nu_{\bxpp\!}^{0,\sigma_n}(x)\\
&=& \liminf nm_p (\nu_{\bxpp\!}^{0,\sigma_n})\\
&=& \liminf \sum_{\pi\in\NC(p)}n\prod_{V\in\pi} k_{|V|}(\nu_{\bxpp\!}^{0, \f{\sigma}{n}})\\
&=& \liminf \sum_{\pi\in\NC(p)}n^{1-\#\pi}\prod_{V\in\pi} k_{|V|}(\nu_{\bxpp\!}^{0, \sigma})\\ &<&\infty,
\end{eqnarray*}
which closes the proof. \ \ \ $\square$

Note that to remove the supplementary hypothesis ``$p$ is even or the support of $\sigma$ is minorized or majorized'' in the second part of the previous \teov it would be useful to prove  proposition P: 

P$:=\{$if $\mu,\nu$ are \pro measures respectively concentrated on $(-\infty,0]$, $[0,\infty)$ \st $\mu\bxpp\nu$ admits a $p^{th}$ moment, then each of them does so$\}$.

 Indeed, in this case, supposing that $\nu_{\bxpp\!}^{\gamma,\sigma}$ admits a $p^{th}$ moment would imply (\cite{defconv}) that there exists free selfadjoint operators $X,Y,Z$ affiliated to a $W^*$-\pro space $(\mc{A},\vfi)$ with respective disributions $\nu_{\bxpp\!}^{0, \sigma^-}$, $\nu_{\bxpp\!}^{\gamma, \sigma^c}$, $\nu_{\bxpp\!}^{0, \sigma^+},$ where $\sigma^-$,$\sigma^c$, $\sigma^+$ are as in the second paragraph  of the proof of $1^\circ$) of the previous proposition, \st $X+Y+Z\in L^p(\mc{A},\vfi)$. Then, since $Y$ is bounded (its distribution has compact support), $X+Z=(X+Y+Z)-Y\in L^p(\mc{A},\vfi)$, so,  by P, $\nu_{\bxpp\!}^{0, \sigma^-}$, $\nu_{\bxpp\!}^{0, \sigma^+}$ admit $p^{th}$ moments. As a consequence, by $2^\circ$) of the previous proposition, $\sigma^-$, $\sigma^+$ admit $p^{th}$ moments, and so does $\sigma$. 

In a more general way, few results give a control on the tails of two \pro measures $\mu$, $\nu$ from the tail of their free convolution (like lemma 3 of \cite{feller}, V.6 for classical convolution).

\section*{Appendix: Taylor expansions in $\Delta_{\alpha,\beta}$}
In this section, we prove a proposition used in the proof of theorem \ref{13.08.04.1}. 
\begin{lemappendix}[Taylor formula in $\mc{H}$] 
Consider $f\in \mc{H}$ and an integer $p$.   \\
(i) If for all $i=0,\ldots, p$, the $i^{th}$ derivate $f^{(i)}$ of $f$ admits in zero a limit $a_i\in\C$, then $f$ admits the Taylor expansion 
$$f(z)=\sum_{i=0}^p\f{a_i}{i!}z^i+o(z^p).$$
(ii) Conversely, if $f$ admits the Taylor expansion $f(z)=\sum_{i=0}^p\f{a_i}{i!}z^i+o(z^p)$, then for all $i=0,\ldots, p$, $$\lim_{z\to 0}f^{(i)}(z)=a_i.$$
\end{lemappendix}
\begin{pr}(i) We prove this result by induction on $p$. For $p=0$, it is obvious. 

Suppose that the result has been proved to rank $p$. Suppose that for all $i=0,\ldots, p+1$, $$\lim_{z\to 0}f^{(i)}(z)=a_i.$$
Replacing $f(z)$ by $f(z)-\sum_{i=0}^p\f{a_i}{i!}z^i$, one can suppose that for all $i$, $a_i=0$. Let us prove $f(z)=o(z^{p+1})$. 

Let us first prove that for all $z$ in the domain $\mc{D}_f$ of $f$ \st the segment $[0,z]$ is contained in $\mc{D}_f$, one has $$f(z)=\int_{[0,z]}f'$$(note that such an integral is defined because $f'$ can be continuously extended in zero).
By Cauchy formula, for all positive $\eps$, one has \begin{eqnarray*}f(z)&=&\int_{[\eps z, z]}f'+f(\eps z)\\
&=&\int_{[0, z]}f'+f(\eps z)-\int_{[0, \eps z]}f'.
\end{eqnarray*}The result follows by letting $\eps$ go to zero (and using the fact that $f$ and $f'$ have null limit in zero).

But by induction hypothesis applied to $f'$, one can write $f'(z)=z^p\eta(z)$, with $\eta\in \h$, $\ds \lim_{z\to 0}\eta(z)=0$.
So \begin{eqnarray*}f(z)&=&\int_{[0,z]}f'\\
&=& z\int_0^1 f'(tz)\ud t\\
&=&z^{p+1}\int_0^1 t^p\eta(tz)\ud t\\ &=&o(z^{p+1}).
\end{eqnarray*}
(ii)  This result is also proved by induction on $p$. For $p=0$, it is obvious. 

Suppose the result to be proved to a rank $p$. Suppose $f$ to admit the Taylor expansion \begin{equation}\label{femmefleureouverte}f(z)=\sum_{i=0}^{p+1}\f{a_i}{i!}z^i+ z^{p+1}\upsilon (z),\end{equation}
where $\upsilon\in\h$, $\ds \lim_{z\to 0}\upsilon(z)=0$.
\\
We already have $\ds \lim_{z\to 0}f(z)=a_0$. By induction hypothesis, it remains to prove that  $$f'(z)=\sum_{i=1}^{p+1}i\f{a_i}{i!}z^{i-1}+ o(z^{p}).$$  So, after differentiation of (\ref{femmefleureouverte}), it suffices to prove that $\ds \lim_{z\to 0} z \upsilon'(z)=0$. 

So let us fix $\alpha>0$, $\eps >0$, and consider $\beta>0$ \st $$\ds \Delta_{2\alpha,\beta}\subset \mc{D}_f,\qquad \sup_{\Delta_{2\alpha,\beta}}|\upsilon(z)|< \eps.$$
Then by the Cauchy inequality, for all $z\in \Delta_{2\alpha,\beta}$, $$|\upsilon'(z)|\leq\f{\eps}{d(z, \C-\Delta_{2\alpha,\beta})}.$$ 
For $z\in \Delta_{\alpha,\beta}$ small enough, the distance of $z$ to $\C-\Delta_{2\alpha,\beta}$ is realized by its orthogonal projection on one of the straight lines $\{x=2\alpha y\}$, $\{x=-2\alpha y\}$. So it is easy to see, with a picture, that if one considers $\theta\in(-\arctan (1/(\!2\alpha)),\arctan (1/(\!2\alpha)))$ \st $z=|z|e^{-i\f{\pi}{2}+i\theta}$, one  has   \begin{equation}\label{222.10.5.05}d(z, \C-\Delta_{2\alpha,\beta})=|z|\sin[\arctan (1/\!\alpha)-|\theta|]>|z|\sin[\arctan (1/\!\alpha)-\arctan (1/(\!2\alpha))].\end{equation}
Hence $$|z\upsilon'(z)|\leq\f{\eps}{\sin[\arctan (1/\!\alpha)-\arctan (1/(\!2\alpha))]}, $$ and we have proved $$ \ds\lim_{\substack{z\to 0\\z\in \Delta_{\alpha,\beta}}} z\upsilon'(z)=0.$$
\end{pr}

\begin{rmq}\label{cecile.douce}A consequence of this lemma is that one can differentiate Taylor expansions of functions of $\h$ and take anti-derivatives of Taylor expansions of functions of $\h$ (as long as anti-derivatives have finite limits in zero).\end{rmq}
\begin{propoappendix}\label{invdeDLdansH}
Consider $f\in \h$ with Taylor expansion $$f(z)=\sum_{i=1}^pa_iz^i+o(z^p)$$without constant term and with  leading coefficient $a_1=1$. Suppose moreover that $f$ is a bijection with inverse $f^{-1}\in\h$. Then $f^{-1}$ admits also a Taylor expansion  of order $p$  without constant term and with  leading term $z$. 
\end{propoappendix}

\begin{pr}We prove this result by induction on the positive integer $p$. 

$\bullet$ For $p=1$, it suffices to prove that $\ds\f{f^{-1}(z)}{z}$ tends to $1$ when $z$ goes to zero non tangeantially. 

a) Let us first prove that for all $\alpha >0$ fixed, for all $\beta>0$ small enough,  \begin{equation}\label{8111.10.05.05}\Delta_{2\alpha,2\beta}\subset \mc{D}_f,\quad\Delta_{\alpha,\beta}\subset f(\Delta_{2\alpha,2\beta}),\end{equation}where $\mc{D}_f$ denotes the domain of $f$. Let us fix $\alpha >0$. Consider $\beta_0>0$ \st $$ \Delta_{3\alpha,\beta_0}\subset \mc{D}_f.$$ It suffices to prove that for $\beta >0$ small enough, for all $\eps\in (0,1)$, \begin{equation}\label{1111.10.05.05}\Delta_{\alpha,\beta}-\Delta_{\alpha,\eps\beta}\subset f(\Delta_{2\alpha,2\beta}-\Delta_{2\alpha,\f{\eps}{2}\beta}).\end{equation} 
 We will prove (\ref{1111.10.05.05}) as a consequence of Rouch\'e's lemma. 
It suffices to prove that for $\beta$ small enough, for all $\eps\in (0,1)$, for all $\omega\in  \Delta_{\alpha,\beta}-\Delta_{\alpha,\eps\beta}$, on the   boundary of $\Delta_{2\alpha,2\beta}-\Delta_{2\alpha,\f{\eps}{2}\beta}$, the inequality $|f(z)-z|<|z-\omega|$ holds. Consider $\beta >0,\eps\in (0,1)$, $\omega\in  \Delta_{\alpha,\beta}-\Delta_{\alpha,\eps\beta}$ and $z$ in the boundary of $\Delta_{2\alpha,2\beta}-\Delta_{2\alpha,\f{\eps}{2}\beta}$. Then either  $|\Re z|=2\alpha|\Im z|$, or $|z|=\f{\eps}{2}\beta$, or $|z|=2\beta$. 
In the first case, by the same arguments as in the proof of (\ref{222.10.5.05}), $$|z-\omega|>|z|\sin[\arctan (1/\!\alpha)-\arctan (1/(\!2\alpha))].$$ 
In the second case, $$|z-\omega|>\f{\eps}{2}\beta=|z|.$$
In the third case, $$|z-\omega|>\f{\beta}{2}=\f{|z|}{2}.$$
Therefore, since $$\ds\lim_{\substack{z\to 0\\ z\in \Delta_{3\alpha,\beta_0}}}\f{f(z)-z}{z}=0,$$ for $\beta >0$ small enough, for all $\eps\in (0,1)$, for all $ \omega \in  \Delta_{\alpha,\beta}-\Delta_{\alpha,\eps\beta}$, on the   boundary of $\Delta_{2\alpha,2\beta}-\Delta_{2\alpha,\f{\eps}{2}\beta}$, the inequality $|f(z)-z|<|z-\omega|$ holds.

b) Note that (\ref{8111.10.05.05}) implies  
$\Delta_{\alpha,\beta}\subset \mc{D}_{f^{-1}}$ and $f^{-1}(\Delta_{\alpha,\beta})\subset \Delta_{2\alpha,2\beta},$ 
where $\mc{D}_{f^{-1}}$ denotes the domain of $f^{-1}$. So we have proved that $f^{-1}(z)$ tends to zero non  tangeantially when $z$ goes to zero non tangeantially.
Since $ \f{f^{-1}(z)}{z}=\f{f^{-1}(z)}{f(f^{-1}(z))}$, and since  $\ds \f{z}{f(z)}$ tends to $1$ when $z$ goes to zero non tangeantially, it implies that $\ds \f{f^{-1}(z)}{z}$ tends to $1$ when $z$ goes to zero non tangeantially.

$\bullet$ Suppose that the result has been proved to  rank $p$. Consider $f\in \h$ with Taylor expansion $$f(z)=\sum_{i=1}^{p+1}a_iz^i+o(z^{p+1})$$without constant term and with  leading coefficient $a_1=1$. Since $\ds \lim_{z\to 0}f^{-1}(z)=0$, by remark \ref{cecile.douce}, it suffices to prove that the derivative of $f^{-1}$ admits a Taylor expansion of order $p$ with constant term $1$. This assertion will be proved by the formula $(f^{-1})'=\ff{f'\circ f^{-1}}$ and the following succession of arguments.
\begin{itemize}
\item[.] By induction hypothesis, $f^{-1}$ admits a Taylor expansion of order $p$  without constant term and with  leading term $z$.
\item[.] By remark \ref{cecile.douce}, $f'$ admits  a Taylor expansion of order $p$  with constant term $1$.
\item[.] Since, as explained in the first step of the induction, $f^{-1}(z)$ goes to zero non tangeantially when $z$ goes to zero non tangeantially, we can compose both Taylor expansions, and we obtain a Taylor expansion of $f'\circ f^{-1}$ of order $p$,  with constant term $1$.
\item[.] It is obvious that if one composes on the left a function of $\h$ which admits a Taylor expansion  of order $p$ with constant term $1$ by a function which admits a Taylor expansion of order $p$ in a neighborhood of $1$, one obtains a function of $\h$ which admits a Taylor expansion of order $p$. 
Therefore, $\ff{f'\circ f^{-1}}$  admits a Taylor expansion of order $p$ with constant  term $1$. 
\end{itemize}    
\end{pr}

Florent Benaych-Georges\\ DMA, \'Ecole Normale Sup\'erieure,\\ 45 rue d'Ulm, 75230 Paris Cedex 05, France\\  e-mail : benaych@dma.ens.fr\\
  http://www.dma.ens.fr/$\sim$benaych
 \end{document}